\documentclass[12pt,a4paper]{article}
\pdfoutput=1 
\usepackage[utf8]{inputenc}
\usepackage{amsthm}
\usepackage{amsmath}
\usepackage{bbm}
\usepackage{amsfonts}
\usepackage{amssymb}
\usepackage{mathtools}
\usepackage[left=2.5cm,right=2.5cm,top=2.5cm,bottom=2.5cm]{geometry}
\usepackage{tabto}
\TabPositions{2 cm, 4 cm, 6 cm, 8 cm}
\usepackage{tocloft}
\usepackage{appendix}
\usepackage{enumerate}
\usepackage{hyperref}
\hypersetup{
  colorlinks   = true, 
  urlcolor     = black, 
  linkcolor    = black, 
  citecolor   = black 
}
 
\newtheorem{theorem}{Theorem}[section]

\newtheorem*{lemma*}{Lemma}
\newtheorem*{proposition*}{Proposition}
\newtheorem*{theorem*}{Theorem}
\newtheorem*{corollary*}{Corollary}
\newtheorem*{claim*}{Claim}

\theoremstyle{definition}
\newtheorem{definition}[theorem]{Definition}

\newtheorem*{definition*}{Definition}

\theoremstyle{remark}
\newtheorem{remark}[theorem]{Remark}

\theoremstyle{definition}

\newtheorem{question}[theorem]{Question}

\mathtoolsset{centercolon} 

\begin{document}

\title{GAGTA 2023 Problem Session}
\author{\footnote{Sections attributed to the person that presented the
  question at the conference. Sections \ref{sec:Bradford},
  \ref{sec:Cashen}, \ref{sec:FFF}, \ref{sec:NMB}, and \ref{sec:Petyt}
  written by named Questioner. Other sections reconstructed from notes
  by CC.} Henry Bradford, Christopher Cashen, Francesco Fournier-Facio,\\
  Nicol{\'a}s Matte Bon, and Harry Petyt}
\date{}
\maketitle

The conference `Geometric and Asymptotic Group Theory with
Applications (\href{http://gagta.org}{GAGTA}) 2023: Groups and
Dynamics' took place at the \href{https://www.esi.ac.at/}{Erwin
  Schr\"odinger Institute} on July 17-21.
These are the problems that were proposed during the Problem Session.

\section{Residual finiteness growth (Henry Bradford)}\label{sec:Bradford}
Recall that a group $G$ is \emph{residually finite} if, 
whenever $g \in G$ is a nontrivial element, 
there exists a finite-index normal subgroup $N$ 
of $G$ which does not contain $g$. 
The \emph{depth} $D(g)$ of $g$ in $G$ is defined to 
be the minimal possible index of such a 
normal subgroup $N$ of $G$. 
In \cite{BouRab}, Bou-Rabee initiated the study of the 
\emph{residual finiteness growth} of finitely generated 
groups. If $G$ is generated by the finite set $S$, 
then the \emph{residual finiteness growth function} 
$\mathcal{F}_G ^S : \mathbb{N} \rightarrow \mathbb{N}$ 
of $G$ (with respect to $S$) is defined by: 
\begin{equation*}
\mathcal{F}_G ^S (n) = \max \lbrace D(g) : \lvert g \rvert_S \leq n \rbrace. 
\end{equation*}
The function 
$\mathcal{F}_G ^S$ only depends on the choice 
of the finite generating set $S$ up to 
a natural notion of equivalence of functions. 
As such, we can talk unambiguously 
about groups of, say, ``polynomial'' or ``exponential'' 
residual finiteness growth. 
A summary of the known bounds 
on $\mathcal{F}_G ^S$ for various familiar groups 
$G$ is given in the Introduction to \cite{BoRaCheTim}. 
One may also ask which functions can arise as 
the residual finiteness growth of some group. 
In this direction, Bou-Rabee and Seward \cite{BoRaSewa} 
constructed groups of arbitrarily 
fast residual finiteness growth. 

\begin{theorem} \label{BRSeThm}
For any increasing function 
$f : \mathbb{N} \rightarrow \mathbb{N}$, 
there exists a residually finite group $G$, 
generated by a finite set $S$, 
such that for all $n \in \mathbb{N}$, 
$\mathcal{F}_G ^S (n) \geq f(n)$. 
\end{theorem}

By slightly modifying the construction of \cite{BoRaSewa}, 
the author strengthened Theorem \ref{BRSeThm} 
by proving a complementary upper bound as follows 
(see \cite{Brad} for a precise statement). 

\begin{theorem}
For any increasing function 
$f : \mathbb{N} \rightarrow \mathbb{N}$ 
obeying a mild ``smoothness'' condition 
and satisfying: 
\begin{equation} \label{RFGLowerBd}
f(n) \geq \exp \big(n \log (n)^2 \log\log(n)^{1+\epsilon} \big)
\end{equation} 
for some $\epsilon > 0$, 
there exists a residually finite group $G$, 
generated by a finite set $S$, 
such that $\mathcal{F}_G ^S$ is equivalent to $f$. 
\end{theorem}

What is not yet well-understood is 
the spectrum of possible ``slow'' 
residual finiteness growth functions. 
For instance, the following remains open. 

\begin{question}
Does there exist a finitely generated residually finite 
group $G$ of \emph{intermediate} residual finiteness growth (that is, 
for which $\mathcal{F}_G ^S$ grows faster than every 
polynomial but slower than any exponential function)? 
\end{question}

One may also ask about groups with residual 
finiteness growth bounded above by a polynomial. 
Finitely generated virtually nilpotent groups all have residual finiteness growth bounded above 
by functions of the form $C \log (n) ^C$. 
The only non-virtually nilpotent groups for which $\mathcal{F}_G ^S$ is known up to equivalence, 
and is at most polynomial, 
are arithmetic groups in some higher-rank Lie group $\mathbf{G}$, and analogues of these over function fields. 
In these cases, the growth is of the form $n^D$, where $D$ is the dimension of the ambient Lie group. 
One may ask whether these are the only 
possible growth types which can arise. 

\begin{question}
Suppose that $G$ is a residually finite group 
generated by a finite set $S$, 
and suppose $G$ is not virtually nilpotent. 
Suppose there exists $C>0$ such that 
$\mathcal{F}_G ^S (n) \leq C n^C$. 
Must there exist $\alpha \in (0,\infty)$ such that 
$\mathcal{F}_G ^S (n)$ is equivalent to $n^{\alpha}$? 
Must $\alpha \in \mathbb{N}$? 
\end{question}

The answer to the last part is probably negative, even for finite-rank nonabelian free groups: 
conjecturally, these have residual finiteness growth of the order of $n^{3/2}$ 
(see \cite{BradThom} Theorem 1.2 and the discussion following). 
We refer the reader to the survey article \cite{DerFerPen} 
for many more open problems about residual finiteness growth and related invariants.

\hyphenation{Chris-topher}
\section{Geometric v.\ random walk boundaries (Christopher Cashen)}\label{sec:Cashen}
Many different boundaries (Gromov, Bowditch, Busemann, visual, Tits, Roller,
Poisson, Martin, Morse, contracting, sublinearly Morse,
quasi-redircting,\dots) have been mentioned at this conference.
The broad questions is how these boundaries are related to one
another.
Many such relationships are known, especially in the hyperbolic
case. For example if $G$ is a hyperbolic group then its Gromov, Morse,
contracting, and sublinearly Morse boundaries all agree, and are
homeomorphic to the Martin boundary associated to a simple random
walk on a Cayley graph of $G$ \cite{Anc88}.
What happens in the non-hyperbolic case?

Recall that the \emph{Gromov boundary} of a hyperbolic space can be
defined as the set of Hausdorff equivalence classes of geodesic ray
emanating from a chosen basepoint, topologized with the `visual
topology', in which two points are close if their representative rays
closely fellow travel for a long time.
Quasi-isometries between hyperbolic spaces induce homeomorphisms of
their Gromov boundaries, so it makes sense to talk about the Gromov
boundary of a hyperbolic group: all of its Cayley graphs are
quasi-isometric to one another, so all of their Gromov boundaries are
homeomorphic to one another.

Quasi-isometry invariance of the boundary does not extend beyond the
hyperbolic setting, even for the case of CAT(0) groups
\cite{CroKle00}.
Charney-Sultan and Cordes \cite{ChaSul15,Cor17} showed that
quasi-isometry invariance can be recovered by restricting to
hyperbolic-like directions:
Recall that a subset of a metric space is called \emph{Morse} if it is
quasi-convex with respect to quasi-geodesics; that is, a
quasi-geodesic segment with endpoints on the set must stay close to
the set, where `close' depends only on the quasi-geodesic constants.
Equivalence classes of Morse geodesic rays at a basepoint make up the
\emph{Morse boundary} as a set.
The direct limit topology used by Charney-Sultan is quite fine. There
is also a coarser topology that is more reminiscent of the visual
topology and is still quasi-isometry invariant \cite{CasMac19}.
With either of these choices of topology, it is easy to see that the
Morse boundary of an infinite, finitely generated group is
homeomorphic to its Gromov boundary when the group is hyperbolic.
When the group is non-hyperbolic the Morse boundary, with either
topology, is either empty or non-compact, and is too small to be a
good random walk boundary: trajectories of the random walk do not
typically stay close to Morse directions.
Rafi-Qing and Tiozzo \cite{QinRaf22,QinRafTio20} enlarge the Morse
boundary to the \emph{sublinearly Morse boundary} by including
geodesic rays whose Morse-ness is allowed to decay sublinearly in
terms of distance from the basepoint.
It is still quasi-isometry invariant, but (when it is non-empty) is large enough to be a
topological model for the Poisson boundary of a simple random walk on
a Cayley graph of the group, in the sense that it has full hitting
measure.

The \emph{Martin boundary} of a simple random walk on a Cayley graph
is another topological model for the Poisson boundary.
It is obtained by embedding the group into the space of normalized positive
superharmonic functions via Martin kernels, and then taking the
frontier of the image.
It is homeomorphic to the Gromov boundary if the group is hyperbolic,
but is mysterious in most other cases. 

\begin{question}\label{question:slMorse_embed_Martin}
Does the sublinearly Morse boundary of a finitely generated group
embed into the Martin boundary?
\end{question}

The answer is `yes' if the group is hyperbolic; it is the non-hyperbolic case that is
interesting.
The map should not be expected to be surjective in the non-hyperbolic case, as flat directions
are missing from the sublinear Morse boundary.

There is a continuous injection from the Morse boundary (with the
direct limit topology) into the Martin boundary \cite{CorDusGek22},
but in the interesting cases the Morse boundary is not compact, so
this map is not necessarily (and in general is not) an embedding.

The topology used in \cite{CasMac19, QinRaf22,QinRafTio20} is in some
sense ad hoc to make the quasi-isometry invariance proof from the
hyperbolic case go through.
A positive answer to this question would say that actually this
topology does coincide with a `naturally occurring' topology. 

A warm-up for Question~\ref{question:slMorse_embed_Martin} would be
the case of toral relatively hyperbolic groups.
In this case the expected answer is `yes', based on \cite{GekGerPot21,DusGekGer20}.

\section{Stability (Francesco Fournier-Facio)}\label{sec:FFF}
A rich and diverse area in modern group theory is broadly known as the study of \emph{stability problems}. Let $\Gamma$ be a countable discrete group and let $\mathcal{G}$ be a family of groups equipped with bi-invariant metrics. The corresponding stability problem asks whether every \emph{almost homomorphism} $\varphi : \Gamma \to G \in \mathcal{G}$ is \emph{close} to a true homomorphism. This can be formalized in two different ways, depending whether one measures distances locally or globally, leading to the notions of \emph{pointwise $\mathcal{G}$-stability} and \emph{uniform $\mathcal{G}$-stability}. \\

The most classical version of this problem is the case when $\mathcal{G} = \mathcal{U}_{op}$ is the family of finite-dimensional unitary groups endowed with the operator norm, commonly known as \emph{Ulam stability} after \cite{ulam}. In this setting, Kazhdan proved that all amenable groups are uniformly stable \cite{kazhdan}, and Burger--Ozawa--Thom proved that certain high-rank lattices are uniformly stable \cite{BOT}; this was recently widely generalized by Glebsky--Lubotzky--Monod--Rangarajan \cite{asymptotic}.

The two most well-studied stability problems nowadays are pointwise stability with respect to the family $\mathcal{S}$ of finite symmetric groups with the normalized Hamming distance, and the family $\mathcal{U}_{HS}$ of finite-dimensional unitary groups with the normalized Hilbert--Schmidt norm. Pointwise $\mathcal{S}$-stability was first studied by Arzhantseva--Paunescu, who proved that $\mathbb{Z}^2$ is pointwise stable \cite{arzhantsevapaunescu}; a standout result is a strong relation between this problem and invariant random subgroups \cite{stability:IRS}. Pointwise $\mathcal{U}_{HS}$-stability is mostly studied from an operator algebraic point of view, due to the strong relation with characters on groups \cite{hadwinshulman1, hadwinshulman2}. \\

In \cite{ultrametric}, the study of stability problems with respect to \emph{ultrametric} families $\mathcal{G}$ was initiated. The ultrametric assumption makes this problem very peculiar, in particular pointwise and uniform stability are equivalent for finitely presented groups \cite[Theorem 1.4]{ultrametric}. The most interesting case is a non-Archimedean version of the various stability problems with respect to finite-dimensional unitary groups, namely the family $\mathcal{G}_p := \{ (\mathrm{GL}_n(\mathbb{Z}_p), d_p) : n \geq 1 \}$, where
\[d_p(A, B) \coloneqq \max\limits_{i, j} |A_{ij} - B_{ij}|_p.\]
This can actually be generalized from the $p$-adic integers to the ring of integers of a non-Archimedean local field, and the theory can be carried through in this generality, with some distinction between the case of positive characteristic and characteristic $0$. While the paper contains several examples \cite[Theorems 1.8 and 1.9]{ultrametric} and a few non-examples \cite[Proposition 1.11, Theorem 1.12]{ultrametric}, these are all quite specific, and some more natural case studies are missing. In particular:

\begin{question}[{\cite[Question 10.6]{ultrametric}}]
\label{q:stability1}
    Is $\mathbb{Z}^2$ $\mathcal{G}_p$-stable?
\end{question}

Note that the ultrametric nature of the setting implies that we do not need to specify whether the stability problem is pointwise or uniform. The stability of $\mathbb{Z}^2$ can be rephrased as the classical problem \emph{are almost-commuting matrices close to commuting matrices?}, which has been studied for years in many different settings \cite{rosenthal:commuting, halmos:commuting, voiculescu:commuting, lin:commuting, FR:lin:commuting, glebsky:commuting}. In this light, and spelling out the definition of the metric $d_p$, Question \ref{q:stability1} can be rephrased as follows:

\begin{question}[{\cite[Lemma 10.7]{ultrametric}}]
    Is the following true? For all $k \geq 1$ there exists $m \geq 1$ such that if $A, B \in \mathrm{GL}_n(\mathbb{Z}_p)$ satisfy $AB \equiv BA \mod p^m$, then there exist $A', B' \in \mathrm{GL}_n(\mathbb{Z}_p)$ such that $A'B' = B'A'$ and $A \equiv A' \mod p^k, B \equiv B' \mod p^k$.
\end{question}

Question \ref{q:stability1} is open also when $\mathbb{Z}_p$ is replaced by the ring of integers of a non-Archimedean local field of positive characteristic, for instance $\mathbb{F}_p[[X]]$. In fact, that case is even more mysterious, since the stability of many finite groups is also open \cite[Section 7.4]{ultrametric}, unlike in characteristic $0$ where all finite groups are known to be stable \cite[Theorem 1.9]{ultrametric}.

\section{Cogrowth (Ilya Gehktman)}\label{sec:Gehktman}
Let $X$ be a proper geodesic Gromov hyperbolic space.
Let $\Gamma\leq \mathrm{Isom}(X)$ be a discrete cocompact subgroup,
and let $N\trianglelefteq\Gamma$ be an infinite normal subgroup.
Let $\delta(N)$ be the exponential growth rate of $N$.
In the special case that $G$ is a free group and $X$ is a Cayley tree, the
quantity $\delta(N)/\delta(G)$ is called the \emph{cogrowth} of the
quotient group $G/N$.
In the more general case it is known \cite{Gri80,Sul87,MatYabJae20,ArzCas20} that $\delta(N)>\frac{1}{2}\delta(\Gamma)$.

\begin{remark}
The question of whether this bound was true in the case that $\Gamma$
is hyperbolic and $X$ is one of its Cayley graph was posed in 1997
\cite{GriHar97} and was open for around 20 years until \cite{MatYabJae20}.
\end{remark}

\begin{remark}
The upper bound, $\delta(N)\leq \delta(\Gamma)$, is also interesting:
equivalence is related to amenability of $\Gamma/N$. \cite{Bro85,Coh82,CouDalSam18,DouSha16,Gri80,RobTap13,Rob05}
\end{remark}

\begin{question}
When is this bound sharp? Do there exist $N_i\trianglelefteq \Gamma$
with $\frac{\delta(N_i)}{\delta(\Gamma)}\to \frac{1}{2}$?
\end{question}
There are partial answers when $X=\mathbb{H}^2$ \cite{BonMatTay12} or $X=\mathbb{F}_2$ \cite{Gri80}
\begin{question}
If $\Gamma$ is a lattice in $\mathrm{Isom}(\mathbb{H}^3)$, does there
exist a sequence  $N_i\trianglelefteq \Gamma$
with $\frac{\delta(N_i)}{\delta(\Gamma)}\to \frac{1}{2}$?
\end{question}
A weaker question:
\begin{question}
Does there exist a sequence $N_i\leq \mathrm{Isom}(\mathbb{H}^3)$ of
discrete subgroups with each $\mathbb{H}^3/N_i$ having finite
injectivity radius and $\delta(N_i)\searrow 1$?
\end{question}
A theorem of Gehktman says one cannot do better
than 1 in the previous question, but it is not known if 1 is realizable.

\section{Embedding left orderable groups (Yash Lodha)}\label{sec:Lodha}
Osin \cite{Osin2010} showed that every countable torsion-free group
embeds into a finitely generated group with precisely two conjugacy
classes.
The question is whether `torsion-free' can be replaced by left
orderable:
\begin{question}
Does every countable left orderable group embed in a finitely
generated left orderable group having precisely two conjugacy classes?
Do there even exist such examples?
\end{question}

A related question of Navas asks:
\begin{question}
Do there exist finitely generated left orderable groups in which each
nontrivial element is distorted?
\end{question}

\section{Schreier growth gap (Nicol{\'a}s Matte Bon)}\label{sec:NMB}
(These questions are motivated from joint work with Adrien Le Boudec.)
Let $G$ be a group generated by a finite set $S$.
Suppose $X$ is a set on which $G$ acts faithfully.
Define $\mathrm{Vol}_{G,X}(n):=\max_{x\in X}|\mathcal{B}_{G,S}(n) x|$, where $\mathcal{B}_{G,S}(n)$ is the ball of radius $n$ in $G$.

Note that one may also consider the Schreier graph of the action, but
it may not be connected since the action is not assumed transitive. The function $\mathrm{Vol}_{G,X}(n)$ measures the volume of the larger ball in the graph. 

Observe that $\mathrm{Vol}_{G,X}(n)$ is bounded above by
$|\mathcal{B}_{G,S}(n)|$ and below by $n$ (provided $G$ is infinite).

Many groups admit faithful actions such that the growth of
$\mathrm{Vol}_{G,X}(n)$ is asymptotically linear. For example it is easy to construct such actions for non-abelian free groups. Furthermore all right-angled 
Artin groups (RAAGs) admit such actions, as observed in \cite{Salo}; as a consequence various groups that are known to embed into RAAGs also do. More examples include Grigorchuk groups, see \cite{Bar-Gri}, topological full groups of subshifts, etc.

We are interested in obstructions, for a finitely generated group $G$, to admit faithful actions such that the function $\mathrm{Vol}_{G,X}(n)$ has low growth. 

\begin{definition}
We say that $G$ \emph{has Schreier growth gap} $f$ if for every
faithful action of $G$ on a set $X$ we have
$\mathrm{Vol}_{G,X}(n)\succeq f(n)$.
The gap is said to be \emph{nontrivial} if $f$ is superlinear.
\end{definition}

There are a few examples. It is well known and easy to see that property (T) implies an exponential Schreier growth gap,  as a consequence of the much stronger fact that the quasi-regular representation in $\ell^2(X)$ has a spectral gap for every infinite $G$-set $X$. (In fact it is enough to have the weaker property FM that every $G$-action which has an invariant mean has a finite orbit \cite{Cor-FM}). 
Another observation which is not difficult is that the existence of a cyclically distorted element in $G$ implies a non-trivial Schreier growth gap. More examples arise among groups of dynamical origin \cite{MB-graph-germs, LB-MB-comm-lemm} and among finitely generated solvable groups \cite{LB-MB-solv}.

\begin{question}
Find some interesting examples among hyperbolic groups. In particular, is there a hyperbolic
group which satisfy a non-trivial Schreier growth gap, but   does not have any infinite subgroup with Property FM?
\end{question}

Note that even though the non-existence of subgroups with Property FM might be difficult to establish, it would still be interesting to find examples among hyperbolic groups without relying explicitly on Property (T) or FM in any way: the question above should be understood in this spirit. 

\begin{question} Are there hyperbolic groups which satisfy a nontrivial Schreier
growth gap, but not an exponential one?
\end{question}

\begin{question}
Does the mapping class group of a closed surface have nontrivial Schreier growth
gap?
\end{question}

\section{$\ell^p$ models for hyperbolic groups (Harry Petyt)}\label{sec:Petyt}
It is a well-known open problem whether all hyperbolic groups are CAT(0). This can be thought of as asking whether hyperbolic groups act geometrically on an $``\ell^2$-type'' space.

A theorem of Lang \cite{lang:injective} states that every hyperbolic group $G$ acts geometrically on an injective metric space $E(G)$. More precisely, $E(G)$ is a locally finite simplicial complex whose simplices are isometric to simplices in $(\mathbf R^n,\ell^\infty)$.

Since we have a geometric action on an ``$\ell^\infty$-type'' space and do not know whether there is a geometric action on an ``$\ell^2$-type'' space, it is natural to weaken the question to ask the whether there is a geometric action on some ``$\ell^p$-type'' space with $p\in[2,\infty)$. There is precedent for this in the world of Banach spaces. Indeed, infinite groups with property~(T) (which includes many hyperbolic groups) cannot admit proper affine actions on Hilbert spaces, but it is a theorem of Yu \cite{yu:hyperbolic} that for each hyperbolic group $G$ there is some $p\in[2,\infty)$ such that $G$ admits a proper affine action on an $\ell^p$--space.

In CAT(0) spaces, each pair of points is joined by a unique geodesic, but in injective spaces geodesics are not unique. A \emph{bicombing} on a geodesic metric space $X$ is a choice of geodesic between each pair of points. That is, for each $x,y\in X$ we have a specified geodesic $\sigma_{xy}:[0,1]\to X$. We say that a bicombing is \emph{convex} if the function $t\mapsto d(\sigma_{xy}(t),\sigma_{x'y'}(t))$ is convex for all $x,x',y,y'\in X$. Every injective space has a natural bicombing with a slightly weaker property than convexity \cite{lang:injective}, and Descombes--Lang showed that if $G$ is a hyperbolic group then the above space $E(G)$ actually has a unique convex geodesic bicombing \cite{descombeslang:convex}.

The strongest bicombing-type property one can have is for the space to be uniquely geodesic and the unique bicombing is convex. This is known as \emph{Busemann convexity}. CAT(0) spaces are Busemann convex, and in fact being CAT(0) can be characterised as being Busemann convex and satisfying a certain four-point inequality \cite{foertschlytchakschroeder:nonpositive}.

Let us say that an \emph{$\ell^p$ complex} is a simplicial complex that is metrised such that each cell is isometric to a simplex of $(\mathbf R^n,\ell^p)$.

\begin{question}
Let $G$ be a hyperbolic group. Does there exist a number $p\in[2,\infty)$ such that there is an $\ell^p$ complex on which $G$ acts geometrically and which 
\begin{itemize}
\item   has a convex bicombing?
\item   is Busemann convex?
\end{itemize}
\end{question}

Simple sufficient conditions and local-to-global type statements for these properties can be found in \cite{alexanderbishop:hadamard,descombeslang:convex,miesch:cartan,haettel:link,haettelhodapetyt:lp}.

\section{Characterizing hyperbolic groups (Davide Spriano)}\label{sec:Spriano}
One might wonder: if $G$ acts by isometries on a hyperbolic space $X$
such that for every infinite order element $g\in G$ the orbit $\langle
g\rangle.x$ is undistorted, must $G$ actually be a hyperbolic
group?
The answer to this question is `no', because there exist distorted,
torsion-free non-hyperbolic subgroups of hyperbolic groups.

\begin{question}
Are there conditions that imply a positive answer?
\end{question}
Obviously making the action geometric would suffice; the goal is to
give up geometricity, but the replace it by sufficiently strong
conditions on cyclic subgroups.

\bibliographystyle{hypersshort}
\bibliography{ref}

\end{document}